\documentclass[a4paper,oneside,12pt]{article}%
\usepackage{amsmath}
\usepackage{amsfonts}
\usepackage{amssymb}
\usepackage{graphicx,color}%

\makeatletter

\newcommand{\Rmnum}[1]{\expandafter\@slowromancap\romannumeral #1@}
\makeatother

\begin{document}

\title{Sublinear expectation linear regression}
\author{
Lu Lin\footnote{The corresponding
author. Email: linlu@sdu.edu.cn. The research was
supported by NNSF projects (11171188, 11071145, 11221061, and 11231005) and the 111 project (B12023) of China,  NSF and SRRF projects (ZR2010AZ001
and BS2011SF006) of Shandong Province of China and K C Wong-HKBU
Fellowship Programme for Mainland China Scholars 2010-11.}, Yufeng Shi, Xin Wang and Shuzhen Yang
\\ Shandong University Qilu Securities Institute for Financial Studies \\and School of Mathematics, Shandong University, Jinan, China}
\date{}
\maketitle

\vspace{-4ex}

\begin{abstract} \baselineskip=16pt Nonlinear expectation, including sublinear expectation as its special case, is a new and original framework of probability theory and has potential applications in some scientific fields, especially in finance risk measure and management. Under the nonlinear expectation framework, however, the related statistical models and statistical inferences have not yet been well established. The goal of this paper is to construct the sublinear expectation regression and investigate its statistical inference. First, a sublinear expectation linear regression is defined and its identifiability is given. Then, based on the representation theorem of sublinear expectation and the newly defined model, several parameter estimations and model predictions are suggested, the asymptotic normality of estimations and the mini-max property of predictions are obtained. Furthermore, new methods are developed to realize variable selection for high-dimensional model. Finally,
simulation studies and a real-life example are carried out to illustrate the new models and methodologies.
All notions and methodologies developed are essentially different from classical ones and can be thought of as a foundation for general nonlinear expectation statistics.

{\it Key words:} Sublinear expectation, $G$-normal, linear regression, model uncertainty, parameter estimation, variable selection.


{\it Running head:} Sublinear expectation regression.

\end{abstract}

\newpage
\baselineskip=21pt

\setcounter{equation}{0}
\section{Introduction}

Among all the assumption conditions imposed to classical statistical models, the most vital one may be that the model under study has a certain probability distribution that may or may not be known. The classical linear expectation and determinant statistics are built on the distribution-certainty or model-certainty. The distribution-certainty, however, is not always the case in practice, such as risk measure and super-hedging in finance. For related references see, e.g., El Karoui, Peng and Quenez (1997), Artzner, Delbaen, Eber and Heath (1999), Chen and Epstein (2002), F\"ollmer and Schied (2004). We also studied a relevant practical problem. It is known that in a financial market, non-performing loan (NPL) is always an important object to be monitored. The NPL ratio is of course related to some economic indicators such as loan-deposit ratio and capital adequacy ratio. We have used an indicator set and the corresponding data published in Vendors Database of China (2000-2010) to establish a regression relationship between the NPL ratio and the indicators in the set. It has been discovered that the regression error has a mean-uncertainty, meanly, the error mean is distributed in an interval $[-0.1833,0.1747]$. We will discuss the issue in detail in Section 5.

Without distribution-certainty, the resulting expectation is nonlinear usually. The earlier works on nonlinear expectation may ascend to Huber (1981) in the sense of robust statistics or ascend to Walley (1991) in the sense of imprecise probabilities. In the recent decades, the theory and methodology of nonlinear expectation have been well developed and received much attention in some application fields such as finance risk measure and control. A typical example of the nonlinear expectation, called $g$-expectation (small g), was introduced in Peng (1997) in the framework of backward stochastic differential equations. As a further development, $G$-expectation (big g) and its related versions are proposed by Peng (2006). Under the nonlinear expectation framework, the most common distribution is the so-called $G$-normal distribution, which was first introduced in Peng (2006). Furthermore, as a theoretical basis of the nonlinear expectation, the law of large numbers as well as the central limit theorem were also established by Peng (2008 and 2009). Also, from different points of view, many authors studied nonlinear expectation, see, e.g., Denis and Martini (2006), Denis et al. (2011), Soner et al. (2011a, 2011b, 2012 and 2013). Other references include Chen and Peng (2000), Briand et al. (2000), Coquet et al. (2002),  Gao (2009), Li and Peng (2011), Peng (1999, 2004, 2005 and 2009), Rosazza (2006), Song (2012), and Xu and Zhang (2009), among many others.

Contrary to the fast development of the nonlinear expectation in probability theory, little attention was paid to the related statistical models and statistical inferences to the best of our knowledge.
Although the earlier work of Huber (1981) refers initially to a upper-lower expectation, a special nonlinear expectation, the main aspects focus on robust statistics and the underlying true model is supposed implicitly to have a certain distribution. Gross error model, for example, contains a certain true distribution in the contaminated distribution set, and based on such a distribution set, the supper-lower expectation can be defined; see, e.g., Strassen (1964) and Huber (1981). In classical statistical frameworks, the heteroscedastic model may be the closest one to the model-uncertainty aforementioned, but it only has variance-uncertainty and the corresponding inference methods do not involve any notion of nonlinear expectation. In nonparametric framework, the model structure is not specified, and in Bayesian framework, the model parameter is random. But the two statistical frameworks are essentially different from the model-uncertainty aforementioned and the corresponding methods are completely unrelated to any nonlinear expectation. In time series models, although the data depend on observation time, strict stationarity or weak stationarity is required to guarantee the certainty of statistical inferences. In a word, under the classical statistical framework, including parameter models, nonparametric models, Bayes models and time series models, the defined expectations are of linearity. Without this linearity, it is essentially difficult or impossible by using classical methods to achieve classical certain conclusions, such as the identifiability of model parameter, estimation consistency, asymptotic normality of the estimation and model selection consistency.

Under model-uncertainty frameworks, the classical statistics methods may no longer be available. The classical maximum likelihood, for example, is nonexistent or can not be uniquely determined due to without a certain likelihood function. Also the classical least squares estimation is invalid because the parameter is defined via linear expectation. Moreover, the classical statistical models such as linear regression models, may not be well-defined as their identifiability depends on mean-certainty; without mean-certainty, the regression function is unidentifiable. Furthermore, it will be verified by simulations in Section 5 that under the situation of model-uncertainty, usual methods may not work and even collapse nearly.
Thus, to achieve the target of statistical inference, it is necessary to develop new statistical frameworks and new statistical methods.

The main contribution of our paper is to establish a framework of sublinear expectation regression for the model that has the distribution-uncertainty.
Based on a sublinear expectation space, a sublinear expectation linear
regression is defined and its identifiability is achieved. Our model is always available for the cases of variance-uncertainty and/or mean-uncertainty. Unlike classical regression, the new model tends to use a large value to predict the response variable and obtains the mini-max prediction risk. It implies that our method is a robust strategy and has potential applications in finance risk measure and management. Based on the representation theorem of sublinear expectation, new parameter estimation methods are suggested and the resulting estimators are asymptotically normal distributed for the case of high-frequency data. Finally, our method is extended to variable selection for high-dimensional regression. It is worth mentioning that under model-uncertainty framework, certainty-statistical inferences are established in this paper, including  parameter-certainty, prediction-certainty and distribution-certainty of parameter estimation. The notions and methodologies developed here are nonclassical and original, and the theoretical framework establishes the foundations for general nonlinear expectation statistics.

The remainder of the paper is organized in the following way. In Section 2, a sublinear regression model is built and its identifiability is obtained. The estimation and prediction methods are suggested in Section 3. Also the asymptotic normality of  estimators and the mini-max property of predictions are established in this section. The method is extended to variable selection for high-dimensional model in Section 4. Simulation studies and a real-life example are carried out in Section 5 to illustrate the new model and methodology. The proofs of the theorems and the definition of the sublinear expectation space are postponed to Appendix.

\setcounter{equation}{0}
\section{Sublinear expectation regression}

In this section we establish a framework of sublinear expectation regression, including modeling, estimation, prediction and asymptotic properties.

\subsection{Model} We consider the following linear regression model:
\begin{eqnarray}\label{(2.1)}
Y =\beta'{\bf x} + \varepsilon,
\end{eqnarray} where $Y$ is a scalar response variable, ${\bf x}=(X_1,\cdots,X_p)'$ is the associated $p$-dimensional covariate having a certain distribution $F_{\bf x}(x)$, and $\beta=(\beta_1,\cdots,\beta_p)'$ is a $p$-dimensional vector of unknown parameters. Furthermore, it is supposed that the error $\varepsilon$ is independent of ${\bf x}$. We need the independence condition only for simplicity. The idea and methodology developed below can be extended to the dependent case, but the notations and algorithm are relatively complex. It is worth pointing out that the essential difference from the classical regression model is that here the error $\varepsilon$ has distribution-uncertainty, which is defined in the following way.

Let $\Omega$ be a given set and $\cal H$ be a linear space of real valued functions defined on $\Omega$. Furthermore, let $\mathbb{E}$ denote a sublinear expectation: $\cal H\rightarrow\mathbb{R}$, satisfying monotonicity, constant preserving, sub-additivity and positive homogeneity; for the details of the definitions see Appendix. The triple $(\Omega,{\cal H},\mathbb{E})$ is then called a sublinear expectation space. In this paper, we assume that the random variable $\varepsilon$ is defined on a sublinear expectation space $(\Omega,{\cal H},\mathbb{E})$.
It can be seen from the definition that the probability distribution of $\varepsilon$ is uncertain. Under this situation, the independence between $\bf x$ and $\varepsilon$ mentioned above is defined in the sublinear expectation space, which is a weak independence (2008 and 2009). For regression analysis, we suppose that $\cal H$ contains linear and quadratic functions, and although the sublinear expectation $\mathbb{E}$ is supposed to be existent, its exact form may be unknown. Thus, a remarkable point of view is that since regression analysis depends mainly on ``expectation", we here only define a sublinear expectation space, instead of the well-accepted linear expectation.

By the representation theorem of sublinear expectation (Peng 2008 and 2009), the sublinear expectation of a function $g(\varepsilon)\in \cal H$ can be expressed as a supremum of linear expectations. Formally, there exists a family of linear expectations $\{E_f:f\in\mathcal{F}\}$ defined on $\cal H$ such that
\begin{eqnarray}\label{(2.1+1)}\mathbb{E}[g(\varepsilon)]=\sup_{f\in\mathcal{F}}
E_f[g(\varepsilon)]\ \mbox{ for } g\in\cal H,\end{eqnarray} and there exists a $f_g\in\mathcal{F}$ such that \begin{eqnarray}\label{(2.1+2)}\mathbb E[g(\varepsilon)]=E_{f_g}[g(\varepsilon)].\end{eqnarray} Write
\begin{eqnarray*}\overline{\mu}=\mathbb{E}[\varepsilon], \ \underline{\mu}=-\mathbb{E}[-\varepsilon], \ \overline\sigma^2=\mathbb{E}[\varepsilon^2], \ \underline\sigma^2=-\mathbb{E}[-\varepsilon^2].\end{eqnarray*}
Then, the intervals $[\underline\mu,\overline\mu]$
and $[\underline\sigma^2,\overline\sigma^2]$ characterize the mean-uncertainty and the variance-uncertainty of $\varepsilon$, respectively.

When ${\bf x}$ is a random variable, for regression modeling, it is necessary to clarify the sublinear expectation $\mathbb{E}[Y|{\bf x}]$ conditional on $\bf x$ since the nonlinear conditional expectation has not yet been defined in the existing literature. Actually, however, there is no obstacle to extend the nonlinear unconditional expectation to the nonlinear conditional expectation. By the representation theorem given above, for instance,
the above $\mathbb{E}[Y|{\bf x}]$ can be defined as
$$\mathbb{E}[Y|{\bf x}]=\sup_{f_{Y|{\bf x}}\in\mathcal{F}_{Y|{\bf x}}}
E_{f_{Y|{\bf x}}}[Y|{\bf x}],$$ where $\{E_{f_{Y|{\bf x}}}:f_{Y|{\bf x}}\in\mathcal{F}_{Y|{\bf x}}\}$ is a family of conditional linear expectations. With this definition, the properties of monotonicity, constant preserving, sub-additivity and positive homogeneity given in Appendix still hold.

Finally, we should note that it was assumed above that the covariate vector $\bf x$ has a certain distribution $F_{\bf x}$ and  the intercept term of model (\ref{(2.1)}) is zero. Here we need the distribution-certainty of $\bf x$ to guarantee that the regression coefficient vector $\beta$ is identifiable; otherwise, when both $\varepsilon$ and $\bf x$ do not have the distribution-certainty, $\beta$ can not be uniquely determined. For details see Remark 2.1 below. The assumption on $\bf x$ and $\varepsilon$ aforementioned is a practical condition. For example, if $Y$ is a measure of a financial risk and $\bf x$ is the set of the corresponding economic indicators, then, usually the goal of regression analysis is to describe the risk measure $Y$ for a given economic indicator set $\bf x$. Therefore, the indicator elements of $\bf x$ could be regarded as of distribution-certainty exactly or approximately. In this case, the model-uncertainty is derived from the unstable financial environments that can be grouped in the model error $\varepsilon$. On the other hand, we need the zero intercept to eliminate the estimation bias; without it, the estimation is inconsistent. For details see Remark 3.2 below.

\subsection{$G$-normal regression}

We first consider the case when the error $\varepsilon$ is supposed to be $G$-normally distributed, namely,
\begin{eqnarray}\label{(2.2)}
\varepsilon \sim \mathcal N=N\left(\{0\}\times[\underline{\sigma}^2,\overline{\sigma}^2]\right).
\end{eqnarray} Under this situation, $\varepsilon$ has a certain zero mean but its variance is uncertain, a special distribution-uncertainty. As was defined by Peng (2006), $\varepsilon$ is called $G$-normally
distributed if it is defined on a sublinear expectation space $(\Omega,{\cal H},{\mathbb E})$ and satisfies that for each $a,b\geq0,$
$$
a\varepsilon+b\,\bar{\varepsilon}\overset{d}{=}\sqrt{a^{2}+b^{2}}\,\varepsilon,
$$
where $\bar{\varepsilon}$ is an independent copy of $\varepsilon$ and ``$\overset{d}{=}$" stands for equal in distribution. For the definition and the representation of $G$-normal distribution see Peng (2006).
It follows from the cash translatability of sublinear expectation given in Appendix that for regression model (\ref{(2.1)}), if $\varepsilon$ is $G$-normally distributed as in (\ref{(2.2)}), then
\begin{eqnarray}\label{(2.4)}
\mathbb{E}[Y|{\bf x}]=\beta'{\bf x}.
\end{eqnarray}
The above relationship (\ref{(2.4)}) could be thought of as a $G$-normal expectation regression because $\mathbb{E}$ is the $G$-normal expectation, a special sublinear expectation. Note that ${\bf x}$ has an identical distribution. Then, we have the following conclusion.

{\bf Proposition 2.1} {\it (1) If $\varepsilon$ is $G$-normally distributed as in (\ref{(2.2)}), then, the $G$-expectation of $Y$ is identifiable in the sense that $\mathbb{E}(Y|{\bf x})$ can be uniquely determined by $\beta'{\bf x}$ as in (\ref{(2.4)}). (2) Besides the condition above, if $E[{\bf x}{\bf x}']$ is a positive definite matrix, then, $\beta$ is identifiable in the sense that $\beta$ can be uniquely determined as
\begin{eqnarray}\label{(2.5)}\beta=(E[{\bf x}{\bf x}'])^{-1}E\{{\bf x}\mathbb{E}[Y|{\bf x}]\},\end{eqnarray} where linear expectation $E$ is taken under the certain distribution $F_{\bf x}(x)$.}

The proof is given in Appendix. For the proposition, we have the following remark.

{\bf Remark 2.1} \begin{itemize}\item[(1)]
The proposition implies that if the error $\varepsilon$ is $G$-normally distributed and $\bf x$ has the distribution-certainty, then $G$-normal regression has both  regression function-certainty and regression coefficient-certainty. The conclusion provides a theoretical basis for regression analysis such as parameter estimation and model prediction.
\item[(2)] From the proof of the proposition we can see that if $\bf x$ does not have the distribution-certainty but only a sublinear expectation is defined for $\bf x$, $\beta$ can not be uniquely determined usually. Without the identifiability of $\beta$, there is no sense in
    modeling regression relationship.
\item[(3)] Here we emphasize the use of $G$-normal regression because a quadratic loss function will be employed below to construct a ``quasi maximum likelihood" estimation; for details see the next section. In fact
the notion proposed here can be directly extended to general mean-certainty sublinear expectation regressions. Specifically, we only assume $\varepsilon$ has the mean-certainty, instead of $G$-normal distribution. Under this situation, model (\ref{(2.4)}) could be regarded as a mean-certainty sublinear expectation regression. With the point of view, the conclusions in Proposition 2.1 still hold.
 \end{itemize}

\subsection{Sublinear expectation regression}

Now we investigate the model in which the error $\varepsilon$ is mean-uncertain and variance-uncertain.
By the cash translatability of the sublinear expectation given in Appendix, we have
\begin{eqnarray}\label{(2.6)}
\mathbb{E}[Y|{\bf x}]=\beta'{\bf x}+\overline\mu.
\end{eqnarray} This model could be thought of as a sublinear expectation regression because $\mathbb{E}$ is a sublinear expectation. By (\ref{(2.6)}) and similar arguments used in Proposition 2.1, we have the following conclusion.

{\bf Proposition 2.2} {\it (1) If $\underline{\mu}<\overline{\mu}$, then, given $\bf x$, the sublinear expectation of $Y$ has a shift $\overline\mu$, more precisely, the sublinear expectation of $Y$ has the framework of (\ref{(2.6)}). (2) Besides the condition above, if $E[{\bf x}{\bf x}']$ is a positive definite matrix, then, $\beta$ is identifiable, more precisely, $\beta$ can be uniquely expressed by
\begin{eqnarray}\label{(2.7)}\beta=(E[{\bf x}{\bf x}'])^{-1}E\{{\bf x}\mathbb{E}[Y|{\bf x}]\}-\overline{\mu}(E[{\bf x}{\bf x}'])^{-1}E[{\bf x}].\end{eqnarray}
Particularly, if $E[{\bf x}]=0$, then
\begin{eqnarray}\label{(2.8)}\beta=(E[{\bf x}{\bf x}'])^{-1}E\{{\bf x}\mathbb{E}[Y|{\bf x}]\}.\end{eqnarray} }

The proof is also presented in Appendix. From the proposition, we have the following findings.

{\bf Remark 2.2}
\begin{itemize} \item[(1)]
The conclusions in the proposition are somewhat surprising because they suggest a nonclassical point of view and provide a methodological development. That is to say, in the face of mean-uncertainty, we can still uniquely determine the parameter vector $\beta$ and then use the mean-shift framework $\beta'{\bf x}+\overline\mu$, instead of $\beta'{\bf x}$, to predict the response variable $Y$. Such a framework reflects the robust feature of sublinear expectation regression. If $Y$ is a measure of the risk of a financial product, then the sublinear expectation regression tends to use a relatively large value to predict risk and moreover, and the increment of risk measure is just the sublinear expectation $\overline\mu$ of the error $\varepsilon$.
\item[(2)]
It is worth mentioning that when the model does not have the mean-certainty, the representation (\ref{(2.7)}) of regression coefficient vector $\beta$ is different from the representation in (\ref{(2.5)}) for the mean-certainty model, in other words, the representation (\ref{(2.5)}) in the mean-certainty model is a special case of the representation (\ref{(2.7)}) with $\overline\mu=0$. This is an essential feature of sublinear expectation regression, i.e., in the mean-uncertainty framework, the regression coefficient vector $\beta$ depends on the nonlinear expectation of error $\varepsilon$. Such a feature is totally different from classical linear expectation regression because in the linear expectation regression framework, the regression coefficient vector $\beta$ has an error-free representation as
$$\beta=(E[{\bf x}{\bf x}'])^{-1}E[{\bf x}Y].$$
\end{itemize}

On the other hand, when $E[{\bf x}]=0$, the parameter representation in (\ref{(2.8)}) is free of $\overline\mu$.
In the following, we mainly focus on the parameter representation in (\ref{(2.8)}) because we will see that without $\overline\mu$, the corresponding estimator of $\beta$ is relatively simple and is asymptotically unbiased.

\setcounter{equation}{0}
\section{Estimation and prediction}

It is supposed in this section that the dimension $p$ of $\beta$ is fixed.
Let $\{(Y_i,{\bf x}_i:i=1,\cdots,N\}$ be a sample from model (\ref{(2.1)}), satisfying
\begin{eqnarray}\label{(3.0)}Y_i=\beta'{\bf x}_i+\varepsilon_i, \ i=1,\cdots,N.\end{eqnarray}
Unlike the classical ones, here $Y_1,\cdots,Y_N$ may have distribution-uncertainty due to the distribution-uncertainty of $\varepsilon_1,\cdots,\varepsilon_N$. Then the corresponding estimation method should be different from the classical ones that only apply to linear expectation regression models.

It seems that we can use (\ref{(2.8)}) to construct the estimator of $\beta$ as it presents a closed expression for $\beta$. However, in the expression, $\mathbb{E}[Y|{\bf x}]$ is a sublinear conditional expectation, like the classical ones, its estimation does involve multivariate nonparametric methods and therefore faces the curse of dimensionality if the dimension $p$ of $\bf x$ is high. To avoid the problem, we now introduce a mini-max method to construct the estimator of $\beta$.

{\bf Case 1.} We first consider the case of $\varepsilon$ having the mean-certainty. Because $Y$ has the sublinear expectation $\beta'{\bf x}$ given ${\bf x}$, theoretically, we should choose $\beta$ so that it can minimize the sublinear expectation loss:
\begin{eqnarray}\label{(3.1)}\mathbb{E}
\left[(Y-\beta'{\bf x})^2\right].\end{eqnarray}
We can easily verify that the above sublinear expectation loss is a convex function function of $\beta$. Thus the optimization problem has a unique global optimal solution. The above is in fact a sublinear expectation least squares. It is worth mentioning that under $G$-normal distribution, we have
that if $\varphi$ is a convex function, then
$$\mathbb{E}[\varphi(\varepsilon)]=\frac{1}{\sqrt{2\pi}}\int_{-\infty}^\infty
\varphi(\overline{\sigma}u)\exp\left\{-\frac{u^2}{2\overline{\sigma}^2}\right\}du,$$
and if $\varphi$ is a concave function, then
$$\mathbb{E}[\varphi(\varepsilon)]=\frac{1}{\sqrt{2\pi}}\int_{-\infty}^\infty
\varphi(\underline{\sigma}u)\exp\left\{-\frac{u^2}{2\underline{\sigma}^2}\right\}du.$$ For details refer to Peng (2006). These imply that under the convex function and concave function spaces, the $G$-normal has density functions $\frac{1}{\sqrt{2\pi}\overline\sigma}
\exp\left\{-\frac{u^2}{2\overline{\sigma}^2}\right\}$ and $\frac{1}{\sqrt{2\pi}\underline\sigma}
\exp\left\{-\frac{u^2}{2\underline{\sigma}^2}\right\}$, respectively. Therefore, the above sublinear expectation least squares could be thought of as a ``quasi maximum likelihood".

To implement the estimation procedure, we need the following assumption:
\begin{itemize}\item[\it C1.] There exists an index decomposition: $I_i,i=1,\cdots,m$, such that when $(ij)\in I_i$, $\varepsilon_{i1},\cdots,\varepsilon_{in_i}$ are independent and have an identical distribution. \end{itemize}
This condition is essentially implied in the conclusion (\ref{(2.1+1)}) of the representation theorem given in Subsection 2.1. Thus, $m$ should be equal to the number of functions in $\mathcal{F}$ if $\mathcal{F}$ only contains finite number of functions; otherwise, $m$ should tend to infinity and in this case, the condition {\it C1} is only an approximation of the true one. We will further weaken {\it C1} and suggest a data-driven decomposition after Theorem 3.1 given below.
From now on we suppose that the numbers of elements in $I_i,i=1,\cdots,m$, are equal, i.e., $n_1=n_2\cdots=n_m=n$, without loss of generality. Because it is assumed that $\varepsilon_{i1},\cdots,\varepsilon_{in}$, are identically distributed, the independence in condition {\it C1} is the same as that in linear expectation framework, instead of the independence in the nonlinear expectation. Here we need  independence only for simplicity. Without the independence assumption, for example, $\varepsilon_{i1},\cdots,\varepsilon_{in}$ are weakly dependent, the conclusions given below still hold; for weakly dependent processes and the properties of estimation see for example Rosenblatt (1956,
1970), Kolmogorov and Rozanov (1960), Bradley and Bryc (1985), and Lu
and Lin (1997). Furthermore, a common decomposition is built according to the observation time order, more precisely,
$\varepsilon_1,\cdots,\varepsilon_N$ are reindexed as $\varepsilon_{ij}=\varepsilon_{(i-1)n+j},i=1,\cdots,m,j=1,\cdots,n$, and then the index sets $I_i$'s are defined as $I_i=\{(ij):j=1,\cdots,n\}$. It is known that in a small time interval, the characteristic of data could be regarded as to be changeless exactly or approximately. Under this point of view, condition {\it C1} is relatively mild. Also we can decompose the index set according to the values of $Y$ in a descending order for example.

Denote by $F_i$ the common distribution function of $\varepsilon_{ij},(ij)\in I_i$. According to the representation theorem of sublinear expectation given in (\ref{(2.1+1)}), sublinear expectation loss (\ref{(3.1)}) can be written as $\max\limits_{1\leq i\leq m}E_{F_i}[(Y-\beta'{\bf x})^2]$ and therefore its empirical version is
\begin{eqnarray}\label{(3.2)}\max_{1\leq i \leq m}\frac 1n\sum_{j=1}^n\left[Y_{ij}-\beta'{\bf x}_{ij}\right]^2.\end{eqnarray}
By minimizing (\ref{(3.2)}), we obtain a mini-max estimator of $\beta$ as
\begin{eqnarray}\label{(3.3)}\hat\beta_G=\arg\min_{\beta\in \mathcal B}\max_{1\leq i \leq m}\frac 1n\sum_{j=1}^n\left[Y_{ij}-\beta'{\bf x}_{ij}\right]^2.\end{eqnarray}
It can be easily verified that $\max\limits_{1\leq i \leq m}\frac 1n\sum_{j=1}^n\left[Y_{ij}-\beta'{\bf x}_{ij}\right]^2$ is a convex function function of $\beta$. Thus the resulting estimator $\hat\beta_G$ is a unique global optimal solution in the above optimization problem. Furthermore, such an estimation procedure can be easily implemented via, for example, genetic algorithm.
Denote $\sigma_i^2=E(\varepsilon^2_{ij})$ for $(ij)\in I_i$ and $\sigma^{2}_{i_*}=\max\limits_{1\leq i\leq m}\sigma_i^2$, and for simplicity, assume that
$$\sigma^{2}_{i_*}>\sigma^{2}_{i} \ \mbox{ for all } i\neq i_*.$$
The mini-max estimator above is asymptotically normally distributed. The following theorem gives the details.

{\bf Theorem 3.1} {\it For the mean-certainty model, if condition {\it C1} holds and $E[{\bf x}{\bf x}']$ is a positive definite matrix and $n\rightarrow\infty$ as $N\rightarrow\infty$, then
$$\sqrt{n}(\hat \beta_G-\beta)\stackrel{ d}\longrightarrow N\left(0,\sigma^{2}_{i_*}(E[{\bf x}{\bf x}'])^{-1}\right)\ \ (N\rightarrow\infty),$$ where $\stackrel{d}\longrightarrow$ stands for convergence in distribution and $ N\left(0,\sigma^{2}_{i_*}(E[{\bf x}{\bf x}'])^{-1}\right)$ is a classical normal distribution. 
}

This theorem establishes the theoretical foundation for further statistical inferences such as constructing confidence intervals and test statistics. From the proof of the theorem given in Appendix we can see that condition {\it C1} can be replaced by the following relatively weak condition:
\begin{itemize}\item[\it C1'.] $\varepsilon_{i_*1},\cdots,\varepsilon_{i_*n}$ are independent and have an identical distribution. \end{itemize} This condition only involves the errors with indexes in $I_{i_*}$. Thus it is relatively common and is implied in (\ref{(2.1+2)}), the second conclusion of the representation theorem.
However, recognizing the fact that the number $n$ of data in each small time slice $I_i$ should be relative large, condition {\it C1} or {\it C1'} only applies to the case of high-frequency data. Moreover, by the two conditions, it is implicitly assumed that the index compositions $I_i,i=1,\cdots,m$, or $I_{i_*}$ are known completely. Under some situations, however, it is difficult or impossible to get such exact compositions in advance. Thus, data-driven decompositions are desired in practice. Now we briefly discuss this issue. By condition {\it C1'}, the proof of Theorem 3.1 and (\ref{(2.1+2)}), the mini-max estimator in (\ref{(3.3)}) can be approximately recasted as
\begin{eqnarray}\label{(3.3+0)}\hat\beta_G=\arg\min_{\beta\in \mathcal B}\frac 1n\sum_{j=1}^n\left[Y_{i_*j}-\beta'{\bf x}_{i_*j}\right]^2.\end{eqnarray} Thus, a simple approach is to identify $I_{i_*}$ or its subset. Let $I_i^0=\{(ij):j=1,\cdots,n^0\}$, $ i=1,\cdots,m^0$, be the initial compositions according to the observation time order for example, where $n^0>p$. Note that in the case of mean-certainty, the common LS estimator $\hat\beta_{LS}$ of $\beta$ is consistent. We then
arrange
$\sum_{j=1}^{n^0}(Y_{ij}-\hat\beta'_{LS}{\bf x}_{ij})^2,i=1,\cdots,m^0$, in the descending order as
$$\sum_{j=1}^{n^0}(Y_{i_1j}-\hat\beta'_{LS}{\bf x}_{i_1j})^2\geq\sum_{j=1}^{n^0}(Y_{i_2j}-\hat\beta'_{LS}{\bf x}_{i_2j})^2\geq\cdots\geq\sum_{j=1}^{n^0}(Y_{i_{m^0}j}-\hat\beta'_{LS}{\bf x}_{i_{m^0}j})^2.$$ From (\ref{(2.1+2)}) we can see that when $n^0$ is relatively small, the index set $I_{i_1}^0=\{(i_1j):j=1,\cdots,n^0\}$ can be chosen as an initial choice of $I_{i_*}$ or a subset of $I_{i_*}$. We then use the data in $I_{i_1}^0$, together with approximate formula (\ref{(3.3+0)}), to build the estimator. Since the data size in $I_{i_1}^0$ may be small, it is necessary to enlarge the initial choice $I_{i_1}^0$. To this end, we consider the following hypothesis test:
$$H_0: \sigma_{2}^2=v_1^2 \ \Leftrightarrow \ H_1: \sigma_{2}^2<v_1^2,$$ where $\sigma_2^2$ is the supposed variance of $\varepsilon_{i_2j}$ for $(i_2j)\in I^0_{i_2}=\{(i_2j):j=1,\cdots,n^0\}$ and $v_1^2=\frac{1}{n^0-p}\sum_{j=1}^{n^0}(Y_{i_1j}-\hat\beta'_{LS}{\bf x}_{i_1j})^2$. Classical methods can used to test the hypothesis $H_0$. If $H_0$ is not rejected, then $I_{i_1}^0\bigcup I_{i_2}^0$ could be chosen as an enlarged choice of $I_{i_*}$. The procedure is repeated until the remainder variances are significantly smaller than $v_1^2$. Also we can use cluster analysis and/or discriminant analysis to achieve this goal.

After the estimator $\hat\beta_G$ is obtained, a natural prediction of $Y$ is
\begin{eqnarray}\label{(3.3+1)}\hat Y=\hat\beta'_G{\bf x}.\end{eqnarray}
If model-uncertainty is ignored and common least squares (LS) method is used to construct the estimator $\hat\beta_{LS}$ of $\beta$, then the LS-based prediction is
\begin{eqnarray}\label{(3.3+2)}\hat Y=\hat\beta'_{LS}{\bf x}.\end{eqnarray} Comparing the two estimators by maximum prediction risk and average prediction risk, we have the following conclusion.

{\bf Theorem 3.2} {\it Under the condition of the mean-certainty, whether the variance-uncertainty exists or not, the following relationships always hold:
\begin{eqnarray*}&&\max_{1\leq i \leq m}\frac 1n\sum_{j=1}^n\left[Y_{ij}-\hat\beta'_G{\bf x}_{ij}\right]^2\leq \max_{1\leq i \leq m}\frac 1n\sum_{j=1}^n\left[Y_{ij}-\hat\beta'_{LS}{\bf x}_{ij}\right]^2,\\&&\frac 1m\sum_{i=1}^m\frac 1n\sum_{j=1}^n\left[Y_{ij}-\hat\beta'_G{\bf x}_{ij}\right]^2\geq \frac 1m\sum_{i=1}^m\frac 1n\sum_{j=1}^n\left[Y_{ij}-\hat\beta'_{LS}{\bf x}_{ij}\right]^2.\end{eqnarray*}}

From the theorem, we have the following finding.

{\bf Remark 3.1} \begin{itemize} \item[] The theorem indicates that sublinear expectation regression is a robust strategy that can reduce maximum prediction risk. Thus, it can be expected that such a regression could be useful for measuring and controlling financial risks.
\end{itemize}

{\bf Case 2.} We now consider the case of $\varepsilon$ having both the mean-uncertainty and the variance-uncertainty. In this case $Y$ has the sublinear expectation $\beta'{\bf x}+\overline\mu$ given ${\bf x}$. Theoretically, we should choose $\beta$ so that it can minimize the sublinear expectation loss:
\begin{eqnarray}\label{(3.4)}\mathbb{E}\left[(Y-\beta'{\bf x}-\overline\mu)^2\right].\end{eqnarray}
However, we cannot directly implement the estimation procedure as $\overline\mu$ is unknown usually. We thus design a profile estimation procedure as follows. Let $\hat\beta$ be an initial estimator of $\beta$, which may be the estimator obtained in Case 1 or by common least squares. When $E[\bf x]=0$, Proposition 2.1 and Proposition 2.2 show that the regression coefficient vectors in Case 1 and Case 2 are equal to each other and thus such an initial estimator is also consistent for Case 2. We then estimate $\overline\mu$ by
$$ \hat{\overline\mu}=\max_{1\leq i\leq m}\frac 1n\sum_{j=1}^n\left[Y_{ij}-\hat\beta'{\bf x}_{ij}\right]$$
and finally estimate $\beta$ by
\begin{eqnarray}\label{(3.5)}\tilde\beta_G=\arg\min_{\beta\in \mathcal B}\max_{1\leq i \leq m}\frac 1n\sum_{j=1}^n\left[Y_{ij}-\beta'{\bf x}_{ij}-\hat{\overline\mu}\right]^2.\end{eqnarray}
Denote $\mu_i=E[\varepsilon_{ij}]$, $\sigma_i^2=E(\varepsilon_{ij}-\mu_i)^2$, $v_{i}^2=\sigma_i^2+(\overline\mu-\mu_i)^2$ and $v^2_{k_*}=\max\{v^2_i:i=1,\cdots,m\}$, and for simplicity, assume $v^2_{k_*}>v^2_i$ for all $i\neq k_*$.
By the same argument as that in Theorem 3.1, we can prove that the estimator $\tilde \beta_G$ is asymptotically normal distributed. The following theorem presents the details.

{\bf Theorem 3.3} {\it For mean-variance-uncertainty, if condition {\it C1} holds, $E[\bf x]=0$ and $E[{\bf x}{\bf x}']$ is a positive definite matrix and $n\rightarrow\infty$ as $N\rightarrow\infty$, then
$$\sqrt{n}(\tilde\beta_G-\beta)\stackrel{\cal D}\longrightarrow N\left(0,\sigma^{2}_{k_*}(E[{\bf x}{\bf x}'])^{-1}\right) \ \ (N\rightarrow\infty).$$ }

For proof of the theorem see Appendix. This theorem establishes a foundation for further statistical inferences and data analyses. Here we also need to check the condition {\it C1}. From the estimation procedure given above, we see that it is asymptotically equivalent to determine two index sets, in which the mean of the error and $\frac{1}{n}\sum_{j=1}^{n}(Y_{k_*j}-\beta'{\bf x}_{k_*j}-\overline\mu)^2$ achieve the maximum values $\overline\mu$ and $v_{k_*}^2$, respectively. The approaches are similar to those used in Case 1 and thus details are omitted here. On the other hand, it is worth pointing out that under the situation of mean-certainty, the condition $E[{\bf x}]=0$ is vital for estimation consistency. The following remark will explain its importance.

{\bf Remark 3.2}
\begin{itemize}\item[] For a model that has the mean-variance-uncertainty, if $E[{\bf x}]\neq0$, then, by the relationship between (\ref{(2.7)}) and (\ref{(2.8)}), we can prove the estimator $\hat{\overline\mu}$ of $\overline\mu$ has an asymptotic bias: $-\overline\mu E[{\bf x}'](E[{\bf x}{\bf x}'])^{-1}E[{\bf x}]$. As a result, if $E[{\bf x}]\neq0$, by the same argument as that used in the proof of Theorem 3.3, it can be verified that the estimator $\tilde\beta_G$ has an asymptotic bias as
$$\mbox{bias}(\tilde\beta_G)=(c\overline\mu-\mu_{k_*})(E[{\bf x}{\bf x}'])^{-1}E[{\bf x}],$$ where $c=1-E[{\bf x}'](E[{\bf x}{\bf x}'])^{-1}E[{\bf x}]$. Furthermore, without $E[{\bf x}]=0$, the bias-correction is essentially difficult because, under the model-uncertainty framework, the law of large numbers can not strictly determine the consistency of sample mean; see Peng (2007 and 2008). On the other hand, the condition $E[{\bf x}]=0$ induces that the intercept term in model (\ref{(2.1)}) should be zero, which implies that if the intercept is nonzero, the estimation bias can not be completely eliminated and thus the estimator is inconsistent.
\end{itemize}

With the estimator, a natural prediction of $Y$ is
\begin{eqnarray}\label{(3.6)}\tilde Y=\tilde\beta_G'{\bf x}+\hat{\overline\mu}.\end{eqnarray}
Similar to the properties in Theorem 3.2, the prediction $\tilde Y$ can obtain the mini-max prediction risk.

{\bf Theorem 3.4} {\it Whether or not the mean-uncertainty and the variance-uncertainty exist, the following relationship always holds:
\begin{eqnarray*}\max_{1\leq i \leq m}\frac 1n\sum_{j=1}^n\left[Y_{ij}-\tilde\beta_G'{\bf x}_{ij}-\hat{\overline\mu}\right]^2\leq \max_{1\leq i \leq m}\frac 1n\sum_{j=1}^n\left[Y_{ij}-\hat\beta'_{LS}{\bf x}_{ij}\right]^2.\end{eqnarray*}}

It shows that our proposal is a robust strategy and is therefore useful for measuring and controlling financial risk. Meanwhile, the simulation study given in Section 5 will verify that when model has mean-variance-uncertainty, the average prediction error of the new method is usually smaller that of the LS method, namely,
$$\frac 1m\sum_{i=1}^m\frac 1n\sum_{j=1}^n\left[Y_{ij}-\tilde\beta'_G{\bf x}_{ij}-\hat{\overline\mu}\right]^2< \frac 1m\sum_{i=1}^m\frac 1n\sum_{j=1}^n\left[Y_{ij}-\hat\beta'_{LS}{\bf x}_{ij}\right]^2.$$ It is because the prediction bias of $\hat\beta'_{LS}{\bf x}$ is between $\underline\mu$ and $\overline\mu$, which is not ignorable, especial for the case of $\underline\mu\,\overline\mu>0$.

\setcounter{equation}{0}
\section{Variable selection}

In this section we focus on the case when the dimension $p=p_N$ tends to infinity as sample size $N$ increases.
Under this situation, model (\ref{(2.1)}) is further supposed to be sparse in the sense that only $d$ components $\beta_{l_k},k=1,\cdots,d$, are
nonzero with $d\ll N$. Without loss of generality, it is assumed that the first $d$ coefficients $\beta_1,\cdots,\beta_d$ are nonzero.

Note that under sublinear expectation framework, the identifiability theory about $\beta$ and $\mathbb{E}[Y|{\bf x}]$ given in Proposition 2.1 and Proposition 2.2 is free of the dimension $p$. Thus, for high-dimensional model, the conclusions in Proposition 2.1 and Proposition 2.2 still hold. With the identifiability, we can investigate variable selection, parameter estimation and model prediction under sublinear expectation framework. For simplicity, we only use the LASSO (Tibshirani (1996) and Zou
(2006)) to achieve our goals. The method developed below can be extended to other penalty methods such as SCAD (Fan and Li (2001) Fan and Peng (2004)) and Dantzig selector (Cand\'es and Tao (2007)).

We first consider the case of $\varepsilon$ having the mean-certainty. The theoretical objective function is defined as
\begin{eqnarray}\label{(4.1)}\mathbb{E}\left[(Y-\beta'{\bf x})^2\right]+\lambda\sum_{k=1}^p|\beta_k|,\end{eqnarray} where $\lambda\geq 0$ is a tuning parameter, which controls the amount of regularization applied to the estimate. Under condition {\it C1}, the empirical version of the above objective function is
\begin{eqnarray}\label{(4.2)}\max_{1\leq i \leq m}\frac 1n\sum_{j=1}^n\left[Y_{ij}-\beta'{\bf x}_{ij}\right]^2+\lambda\sum_{k=1}^p|\beta_k|.\end{eqnarray} By minimizing (\ref{(4.2)}), we can achieve the goals of variable selection and parameter estimation simultaneously. It can be verified easily that the objective function above is a convex function of $\beta$. Then, the global minimum solution exist uniquely. Furthermore, such an optimization procedure can be easily implemented via, for example, genetic algorithm. Denoted by $\hat\beta_G$ the solution of the optimization problem (\ref{(4.2)}). Note that most  components of $\hat\beta_G$ are shrunk to zero by choosing a suitable tuning parameter $\lambda$. Then, the goal of variable selection can be realized. After variable selection and parameter estimation being completed, a natural prediction of $Y$ can be chosen as
\begin{eqnarray}\label{(4.3)}\hat Y=\hat\beta'_G{\bf x}.\end{eqnarray}
Similar to the arguments in Theorem 3.2, our method is a robust strategy because the selected model can reduce the maximum prediction risk. Thus, the selected model by sublinear expectation can be employed to measure and control financial risks.

From the proof of Theorem 3.1, we see that when $n$ is large enough, the term of order $O_p(1/n)$ can be ignored and the objective function above is approximately equal to \begin{eqnarray}\label{(4.4)}\frac 1n\sum_{j=1}^n\left[Y_{i_*j}-\beta'{\bf x}_{i_*j}\right]^2 +\lambda\sum_{k=1}^p|\beta_k|,\end{eqnarray} where $i_*$ is the index of the interval $I_{i_*}$ in which the variance of $\varepsilon$ achieves the maximum value. This representation implies that the properties of variable selection and parameter estimation, such as the selection consistency and the Oracle property of the estimator, are the same as those of the standard LASSO. So it is unnecessary to restudy these theoretical properties under the sublinear expectation framework.
However, this representation shows that the number of data in each small time slice $I_i$ should be relative large. Therefore our method only applies to high-frequency data.

If $\varepsilon$ possesses both the mean-uncertainty and the variance-uncertainty, as was shown in the previous selection, we need the condition $E[{\bf x}]=0$ to guarantee the consistency of estimation. Variable selection and parameter estimation can be obtained by minimizing the following objective function
\begin{eqnarray}\label{(4.5)}\max_{1\leq i \leq m}\frac 1n\sum_{j=1}^n\left[Y_{ij}-\beta'{\bf x}_{ij}-\hat{\overline\mu}\right]^2+\lambda\sum_{k=1}^p|\beta_k|.\end{eqnarray} Here $\hat{\overline\mu}$ is an initial estimator of $\overline\mu$ defined by
$$ \hat{\overline\mu}=\max_{1\leq i\leq m}\frac 1n\sum_{j=1}^n\left[Y_{ij}-\hat\beta'_G{\bf x}_{ij}\right],$$ where $\hat\beta_G$ may be the solution by minimizing (\ref{(4.2)}).
Denote by $\tilde \beta_G$ the corresponding solution. Then a prediction of $Y$ is chosen as
\begin{eqnarray}\label{(4.6)}\tilde Y=\tilde\beta_G'{\bf x}+\hat{\overline\mu}.\end{eqnarray} Also this prediction achieves the mini-max prediction risk and the prediction value tends to be larger.

\setcounter{equation}{0}
\section{Simulation study and real data analysis}

\subsection{Simulation study}

In this section we present several simulation examples
to compare the finite sample performances of the sublinear expectation regression
proposed in this paper with the existing competitors, such as the classical LS regression and the LASSO regression. To
get comprehensive comparisons, we use the mean square error (MSE), maximum prediction error (MPE) and average prediction error (APE), together with scatter plots of the estimation and prediction, to assess the different methods. From the simulations given below, we will get the following findings: (1) The new methods can significantly reduce the MPE under all the situations; (2) When the model has the mean-certainty, the advantages of the new methods over the classical LS methods are not very obvious; (3) For the case of the mean-uncertainty, the predictions of the classical LS methods do not work and even collapse nearly, but the new methods can get a valid prediction because the impact of the mean-uncertainty on the new methods can be successfully eliminated by the use of the sublinear expectation of the error. Thus, our proposals are robust to the uncertainties of mean and variance and particularly, for the case of the mean-uncertainty, the advantages of ours are rather obvious.

{\it Experiment 1.} We first consider the following simple linear model
$$Y=\beta_1X_1+\beta_2X_2+\beta_2X_2+\varepsilon.$$ In the simulation procedure, the regression coefficients are chosen as $\beta_k=1,k=1,2,3$, the observation values of $X_k$ are independent and identically distributed from $N(10,2),k=1,2,3$. We choose $\varepsilon\sim N(\{0\}\times[0,3])$, a $G$-normal distribution with certain zero mean. In this case, the model has the mean-certainty. The following way is used to generate the data of $G$-normal distribution approximately. Firstly, generate variance values $\sigma_i^2,i=1,\cdots,m$, from the uniform distribution $U[0,3]$, and then generate the values $\varepsilon_{ij},j=1,\cdots,n$, of $\varepsilon$ from the common normal distribution $N(0,\sigma_i^2)$.
For $m=10$ and $n=10$, the simulation results are reported in Table 1, in which MSE, MPE and APE denote the mean squared error, maximum prediction error and average prediction respectively; for the definitions of MPE and APE see Proposition 2.2. It is clear by the simulation results that the MSE and APE of common LS estimation $\hat\beta_{LS}$ are significantly smaller than those of the $G$-normal estimation $\hat\beta_G$. Such a result is not surprising because, under the mean-certainty model, the common LS estimation $\hat\beta_{LS}$ is consistent but the construction of the new estimation $\hat\beta_G$ only uses the data in a small time interval (essentially, the number of the data used to construct the estimator $\hat\beta_G$ is only 10).
On the other hand, the MPE by the new one $\hat\beta_G$ is significantly small than that by the LS estimator $\hat\beta_{LS}$, which implies than the new method can reduce the maximum prediction risk and therefore is a robust strategy.

\begin{table}[h]
\caption{Simulation results of estimation and prediction for Experiment 1 with $m=10$ and $n=10$}
\label{tab:1} \vspace{0.3cm} \center
\begin{tabular}{c|ccc||c|c}
  \hline\hline
  MSE & $\beta_1$ & $\beta_2$ & $\beta_2$&MPE&APE\\
    \hline
  $\hat\beta_G$   & 0.0080 & 0.0301 & 0.0315& 6.0259&3.5388 \\
  \hline
 $\hat\beta_{LS}$   & 0.0026 & 0.0045 & 0.0037 &6.6122&2.8584\\
  \hline\hline
\end{tabular}
\end{table}

%
%

The simulation results above indicate that when model has the mean-certainty, the advantages of the new methods over the common LS are not rather obvious. Moreover, the new methods even have the disadvantage of instability. In the following, we will see that when model has the mean-uncertainty, our new methods have rather clear advantages over the LS based methods.

{\it Experiment 2.} We reconsider the linear model
$$Y=\beta_1X_1+\beta_2X_2+\beta_2X_2+\varepsilon,$$ which is the same in form as in Experiment 1. However, here the model has the mean-variance-uncertainty as $\varepsilon\sim N([3,5]\times[0,4])$. The other experiment conditions are designed as $X_k\sim N(0,1),k=1,2,3$, $m=10$ and $n=20$. The values of $\varepsilon$ are generated by the following way. Firstly, the values $\mu_i$ of the mean and the values $\sigma_i^2$ of the variance are generated from the uniform distributions $U[3,5]$ and $U[0,4]$ respectively, and then the values $\varepsilon_{ij},j=1,\cdots,n$, of $\varepsilon$ are generated from the common normal distribution $N(\mu_i,\sigma_i^2)$ for $i=1,\cdots,m$. The simulation results are listed in Table 2. For the MSE of the parameter estimation, the results are similar those in Experiment 1, i.e., the MES of the LS estimation is smaller than that of the new estimation because the new method only uses the data in a small subinterval in principle. However, when the mean-uncertainty and variance-uncertainty appear in the model, both the MPE and the APE of the new one are significantly smaller than those of the LS estimator. Particularly, the prediction by the LS seems to be totally invalid. It indicates that ignoring the model-uncertainty will lead to a serious prediction risk.

\begin{table}[h]
\caption{Simulation results of estimation and prediction for Experiment 2 with $m=10$ and $n=20$}
\label{tab:3} \vspace{0.3cm} \center
\begin{tabular}{c|ccc||c|c}
  \hline\hline
  MSE & $\beta_1$ & $\beta_2$ & $\beta_2$&MPE&APE\\
    \hline
 $\tilde\beta_G$   &0.1258  & 0.2769  &0.2398  &14.4254  & 6.8787 \\
  \hline
 $\hat\beta_{LS}$  & 0.1141 &  0.1891& 0.1879 & 36.0253  & 21.2932\\
  \hline\hline
\end{tabular}
\end{table}

{\it Experiment 3.} In this experiment, we consider the following high-dimensional linear model
$$Y=\sum_{j=1}^p\beta_jX_j+\varepsilon.$$ In the simulation procedure we choose $p=40$,
$\beta_j=1$ for $1\leq j\leq 5$ and $\beta_j=0$ for all $j\geq6$, $X\sim N_{40}(0,I_{40})$, $\varepsilon\sim N(\{0\}\times [1,4])$ and the sample size satisfies $m=10$ and $n=200$. Like the condition in Experiment 1, this model has the mean-certainty. We consider the common LS and $G$-normal estimation, as well as use the common LASSO and the $G$-normal LASSO (G-LASSO) defined in Section 4 to select variables and estimate parameters simultaneously. The tuning parameter $\lambda$ is determined by the CV. Under the above experiment condition, for the common LASSO estimation, the value of $\lambda$ is chosen as $\lambda_{LS}=0.0604$; for the G-LASSO, the value of $\lambda$ is chosen as $\lambda_{G}=0.3377$. The simulation results are reported in Table 3 and Figure 1. In Table 3, GNR, LSR, Lasso-GNR and Lasso-LSR stand for the $G$-normal regression, LS regression, LASSO-$G$-normal regression and LASSO-LS regression, respectively. The simulation results in Table 3 can verify that the new methods can efficiently reduce maximum prediction error. From Figure 1 we have the following findings: (1) The the LS methods are more stable than the new methods; (2) Like the common LASSO, the $G$-normal LASSO can efficiently select the active variables.

\newpage

\begin{table}[h]
\caption{Simulation results of prediction for Experiment 3 with independent covariates}
\label{tab:4} \vspace{0.3cm} \center
\begin{tabular}{c|ccccc}
  \hline\hline
 Models & GNR & LSR & Lasso-GNR&Lasso-LSR\\
    \hline
  MPE & 4.0365 & 4.5443 & 4.0221 & 4.1206 & \\\hline
  APE &2.5334 & 2.0801 & 2.2363 & 2.0573\\
  \hline\hline
\end{tabular}
\end{table}

To further examine the behaviors, here we consider the correlated covariates: $X\sim N_{40}(0,\Sigma)$, where $\Sigma$ is $40\times 40$ matrix with the $(ij)$-element as
$$\Sigma_{ij}=\left\{\begin{array}{ll}1,& \mbox{for } i=1,\vspace{1ex}\\
0.5,&\mbox{for } i\neq j.\end{array}\right.$$ The other experiment conditions are designed as the same as the above. The simulation results are presented in Figure 2. The performances of the figures in Figure 2 are similar to those in Figure 1,
but they are not as stable as before because of the correlation among the covariates.

\begin{figure}
\begin{center}
\begin{tabular}[b]{c}
{\includegraphics[width=12cm,height=5cm]{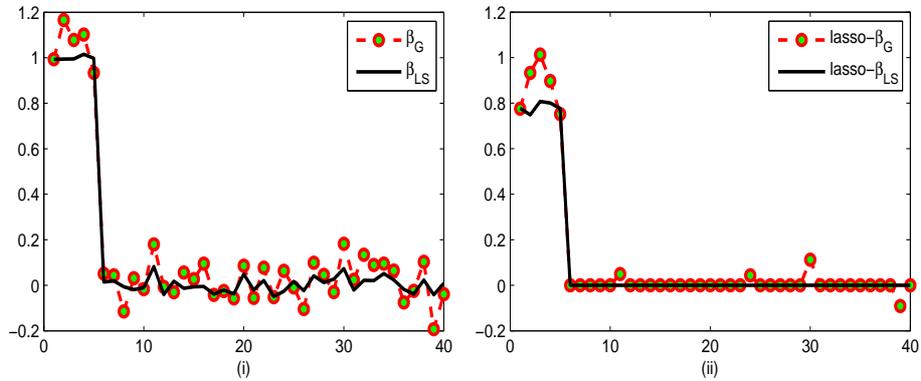}}
\end{tabular}
\caption{The figures of estimation for Experiment 3 with independent covariates.} \label{fig:1}
\end{center}
\end{figure}

\newpage

{\it Experiment 4.}
In this experiment, the model is designed as the same in form as that in Experiment 3, but the model has both the mean-uncertainty and the variance-uncertainty. Formally, the error is distributed as $\varepsilon\sim N([5,10]\times [1,4])$, which has the mean-variance-uncertainty. We first consider the simulations for the GNR and the LSR without use of the LASSO, the results being reported in Figure 3. Figure 3(i) verifies again that the parameter estimation of the LSR is more stable than that of the GNR. On the other hand, Figure 3(ii) provides a clear evidence that with the mean-variance-uncertainty, the LSR has rather large values of the MPE and the APE and therefore the LSR prediction is invalid completely, but the GNR can significantly reduce both the MPE and the APE. These results imply that under the mean-variance-uncertainty framework, ignoring the mean-uncertainty will result in a serious prediction risk, but the new method can efficiently reduce prediction risk by the use of the information of the mean-uncertainty of the error $\varepsilon$.

\begin{figure}
\begin{center}
\begin{tabular}[b]{c}
{\includegraphics[width=12cm,height=5cm]{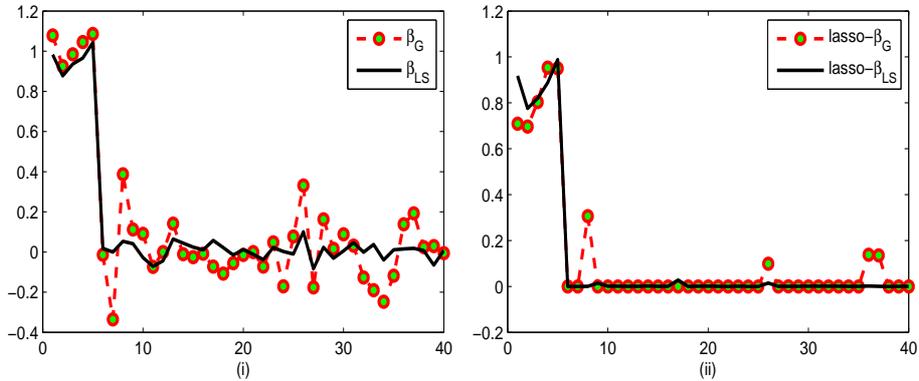}}
\end{tabular}
\caption{The figures of estimation for Experiment 3 with correlated covariates.} \label{fig:2}
\end{center}
\end{figure}

\begin{figure}
\begin{center}
\begin{tabular}[b]{c}
{\includegraphics[width=12cm,height=5cm]{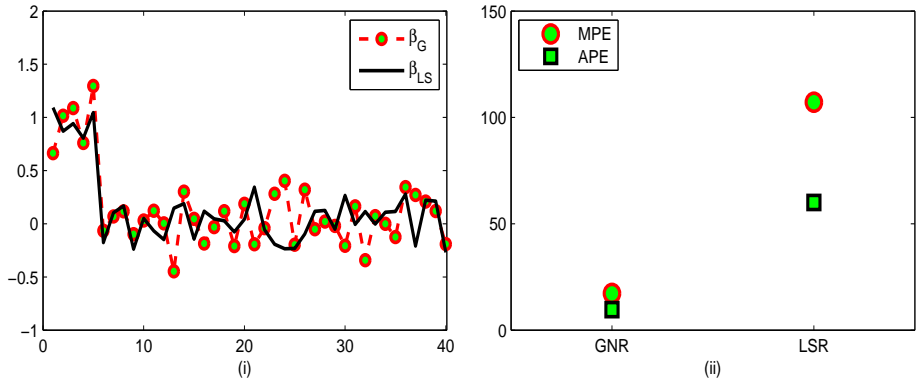}}
\end{tabular}
\caption{The figures of estimation and prediction for Experiment 4.} \label{fig:3}
\end{center}
\end{figure}

Now we consider variable selection and parameter estimation by the LASSO. Under the experiment conditions above, we get $\lambda_{G}=0.6494$ and $\lambda_{LS}=0.4670$ via the CV method. For the Lasso-GNR, the simulation results are given by Figure 4. It shows that the new method can efficiently select active variables and at the same time, the prediction risks are rather small. For the Lasso-LSR, the simulation results are presented in Figure 5. By comparing Figure 5 and Figure 6, we have a clear evidence to show that the new method can obviously reduce the prediction risk, but the LS prediction collapses nearly.

\begin{figure}
\begin{center}
\begin{tabular}[b]{c}
{\includegraphics[width=12cm,height=5cm]{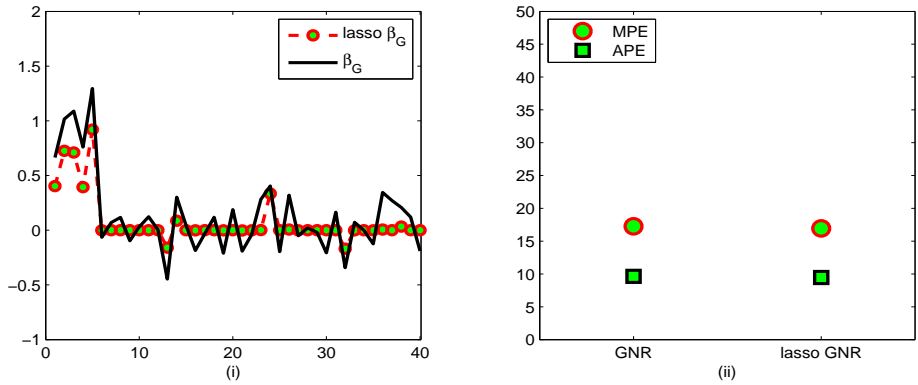}}
\end{tabular}
\caption{The figures of GNR and Lasso-GNR estimation and prediction for Experiment 4.} \label{fig:5}
\end{center}
\end{figure}

\begin{figure}
\begin{center}
\begin{tabular}[b]{c}
{\includegraphics[width=12cm,height=5cm]{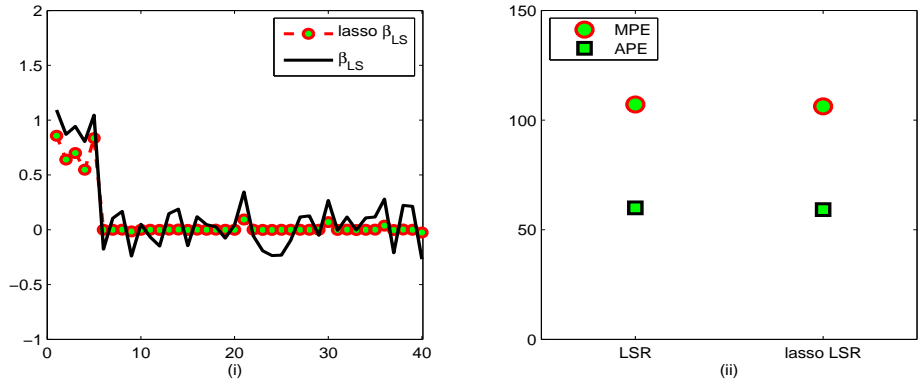}}
\end{tabular}
\caption{The figures of LSR and Lasso-LSR estimation and prediction for Experiment 4.} \label{fig:5}
\end{center}
\end{figure}

\newpage

In short, our methods are robust to the uncertainties of mean and variance. Particularly, for the case of serious mean-uncertainty, the classical methods may collapse, but our new methods can successfully eliminate the impact of mean-uncertainty and construct efficient prediction. The main disadvantage of the new methods is the instability, more precisely, the resulting estimation has relatively large variance since the mini-max estimation only uses the data in a subinterval, essentially.

\subsection{Real data analysis}

Non-performing loan (NPL) is always an important object to be monitored in financial market. To investigate the relationship between the NPL ratio and a set of economic indicators, we use our models, together with the new estimation methods, to fit the real data published in Vendors Database of China (2000-2010). We also compare our fittings with the LS fittings that ignore the distribution-uncertainty. According to the indicator system in Vendors Database, after the indicators with which the data are incomplete are deleted, we choose the following indicators as initial choices: loan-deposit Ratio ($X_1$), capital adequacy ratio ($X_2$), core capital adequacy ratio ($X_3$), liquidity ratio of short-term assets of RMB ($X_4$),
liquidity ratio of short-term assets of foreign currencies ($X_5$), proportion of loans from other banks ($X_6$), proportion of loans to other banks ($X_7$), ten largest customers loan ratio ($X_8$), single biggest customer loan ratio ($X_9$) and NPL provision coverage ($X_{10}$). Because the indicators $X_j$ are percentages, they are transformed to $\tilde X_j=\log\frac{a_j+X_j}{b_j-X_j}$ for some constants $a_j>0$ and $b_j>1$, and then $\tilde X_j$ are centralized so that the centralized versions of $\tilde X_j$ have zero mean. In the following, we still use $X_j$ to denote the transformed and centralized indicators for simplicity. According the observation time order, the data are decomposed into five sets, in which the numbers of valid data are $n_1=26$, $n_2=25$, $n_3=21$, $n_4=20$ and $n_5=31$ respectively.

From the real data analyses given below, we will have the following findings: (1) With model-uncertainty technique, the new methods in most cases have more efficient fitting than the LS does; (2) Particularly, when the technique of mean-variance uncertainty is employed to fit the real data, a more precise fitting can be obtained.

\subsubsection{\bf Case 1 (Mean-certainty model)} We first use a model with mean-certainty to fit the data.

(1) If the variable selection is not taken into account, by our method of variance-uncertainty, we get an empirical model as
\begin{eqnarray*}M_G\mbox{-1}:\ Y=-0.2602X_1+0.1922X_2-0.3953X_3-0.2513X_4+0.0607X_5
\\ -0.1808X_6+0.0727X_7
+0.4314X_8-0.1503X_9-0.5656X_{10}.\end{eqnarray*} With this model, the maximum prediction error and average prediction error have values:
$$MPE(M_G\mbox{-1}) = 1.6009,\ APE(M_G\mbox{-1}) = 0.4632.$$
As a contrastive method, the LS is used to build model, the resulting empirical model has the following form:
\begin{eqnarray*}M_{LS}\mbox{-1}:\ Y=-0.2590X_1+0.1843X_2-0.3972X_3-0.2268X_4+0.0543X_5
\\ -0.2073X_6+0.0914X_7
+0.2734X_8-0.0315X_9-0.5884X_{10}.\end{eqnarray*} The corresponding prediction errors have the following values:
$$MPE(M_{LS}\mbox{-1}) = 1.4396,\ APE(M_{LS}\mbox{-1}) = 0.4323.$$
By comparing the prediction errors, we see that in this case our method has no advantage over the LS fitting. We will analyze the causes in the following studies.

(2) Since some indicators among the ten economic indicators have clear correlation and the number of data is relatively small, the fittings above are inefficient. It is necessary to select variables so that the final model is parsimonious and workable. Now we use the Lasso, together with variance-uncertainty, to build an empirical model, which has the following form:
\begin{eqnarray*}M_G\mbox{-2}:\ Y=-0.1770X_1-0.0111X_2-0.1878X_3\\
-0.0549X_4+0.1397X_8-0.5956X_{10}.\end{eqnarray*} By this treatment, the prediction errors have the following values:
$$MPE(M_G\mbox{-2}) =0.8443,\ APE(M_G\mbox{-2}) = 0.3500.$$
By use of the Lasso, the inactive predictors are removed from the model, the model size is significantly reduced and prediction effectiveness is improved clearly.

If variance-uncertainty is ignored, the Lasso-LS empirical model has the following form:
\begin{eqnarray*}M_{LS}\mbox{-2}:\ Y=-0.0387X_1-0.0269X_2-0.0542X_3\\+0.0722X_8+0.0352X_9
-0.5381X_{10},\end{eqnarray*} and the corresponding prediction error have the following values:
$$MPE(M_{LS}\mbox{-2}) =1.0420,\ APE(M_{LS}\mbox{-2}) = 0.4258.$$
By comparing $M_{G}\mbox{-2}$ and $M_{LS}\mbox{-2}$, we have a clear evidence that our method can reduce prediction errors.

\subsubsection{\bf Case 2 (Mean-variance-uncertainty model)} We can verify that $\underline\mu=-0.1833$ and $\overline\mu=0.1747$. Thus, such a mean-uncertainty is not ignorable. To improve data fitting, both mean-uncertainty and variance-uncertainty are taken into account in the following modeling procedure.

(1) Without use of variable selection, the model with mean-variance-uncertainty has the following empirical expression:
\begin{eqnarray*}\widetilde M_G\mbox{-1}:\ Y=-0.2315X_1+0.1888X_2-0.4765X_3-0.2673X_4+0.0129X_5
-0.2590X_6\\+0.0798X_7
+0.6331X_8-0.3093X_9-0.5374X_{10}+0.1747.\end{eqnarray*}
This model leads to the prediction errors as
$$\widetilde{MPE}(M_G\mbox{-1}) = 0.9182,\ \widetilde{APE}(M_G\mbox{-1}) = 0.3837.$$
Comparing $\widetilde M_G\mbox{-1}$ with both $M_G\mbox{-1}$ and $M_{LS}\mbox{-1}$, the model $\widetilde M_G\mbox{-1}$ has the following two distinctive features: it uses a relatively large value to predict the NPL ratio and the prediction errors are significantly reduced.

(2) By use of the Lasso, the model with mean-variance-uncertainty has the following empirical expression
\begin{eqnarray*}\widetilde M_G\mbox{-2}:\ Y=-0.0389X_1-0.0420X_2-0.1309X_3
-0.5108X_{10}+0.1747.\end{eqnarray*} By this treatment, the prediction errors are reduced to
$$MPE(M_G\mbox{-2}) =0.7305,\ APE(M_G\mbox{-2}) = 0.4362.$$
This model may be the best one among all the models mentioned above because it has both the smallest model size and the smallest MPE.

In short, a flexible model that has mean-variance-uncertainty can relatively precisely fit the real data and is parsimonious and workable.

\setcounter{equation}{0}
\section{Appendix}

\subsection{Definition of sublinear expectation}

Let $\Omega$ be a given set and $\cal H$ be a linear space of real valued functions defined on $\Omega$.
Suppose that $\mathbb{E}:\cal{H}
\rightarrow \mathbb{R}$ satisfies the following properties: for all
$U,V\in \mathcal{H}$,
\begin{itemize}
\item[(i)] Monotonicity: If $U\geq V$ then $\mathbb{E}[U]\geq \mathbb{E}[V]$;
\item[(ii)] Constant preservation: $\mathbb{E}[c]=c$ for any constant $c$;
\item[(iii)] Sub-additivity: $\mathbb{E}[U+V]\leq \mathbb{E}[U]+\mathbb{E}[V]$;
\item[(iv)] Positive homogeneity: $\mathbb{E}[\lambda U]=\lambda
\mathbb{E}[U]$ for each $\lambda \geq0$.
\end{itemize}
Then $(\Omega,\mathcal{H},\mathbb{E})$ is called a sublinear expectation space.

It can be verified that (iii) and (iv) together imply
\begin{itemize}\item[(v)] Convexity:
$$\mathbb{E}[\alpha U+(1-\alpha)V]\leq \alpha\mathbb{E}[U]+(1-\alpha)\mathbb{E}[V] \ \mbox{ for } \alpha\in[0,1].$$\end{itemize}
Furthermore, (ii) and (iii) together lead to
\begin{itemize}\item[(vi)] Cash translatability:
$$\mathbb{E}[U+c]=\mathbb{E}[U]+c \ \mbox{ for any constant } c.$$\end{itemize}

\subsection{Proofs}

{\bf Proof of Proposition 2.1}
We only need to prove the second result. It is clear that (\ref{(2.4)}) yields
$${\bf x}\mathbb{E}[Y|{\bf x}]={\bf x}{\bf x}'\beta.$$ Note that the distribution of $\bf x$ is certain. Consequently,
$$E\{{\bf x}\mathbb{E}[Y|{\bf x}]\}=E[{\bf x}{\bf x}']\beta.$$ This implies the second result of the proposition. $\Box$

{\bf Proof of Proposition 2.2}
We only need to prove the second result. It is obvious that by (\ref{(2.4)}) we have
$${\bf x}\mathbb{E}[Y|{\bf x}]={\bf x}{\bf x}'\beta+\overline{\mu}{\bf x}$$ and consequently
$$E\{{\bf x}\mathbb{E}[Y|{\bf x}]\}=E({\bf x}{\bf x}')\beta+\overline{\mu}E[{\bf x}].$$ This implies the second result of the proposition. $\Box$

{\bf Proof of Theorem 3.1} It follows from {\it C1} that
$$\frac{1}{n}\sum_{j=1}^n\varepsilon_{ij}^2=\sigma^2_i+\delta_n,$$ where $\delta_n$ is of order $O_p(1/n)$ and is free of $\beta$. Consequently,
$$\max_{1\leq i\leq m}\frac{1}{n}\sum_{j=1}^n\varepsilon_{ij}^2=\sigma_{i_*}^2+\delta_n.$$  Denoted by $\beta^0$ the true value of $\beta$. Then
\begin{eqnarray*}&&\max_{1\leq i\leq m}\frac{1}{n}\sum_{j=1}^n\left[Y_{ij}-\beta'{\bf x}_{ij}\right]^2\\&&=\max_{1\leq i\leq m}\frac{1}{n}\sum_{j=1}^n\left[\varepsilon_{ij}-(\beta-\beta^0)'{\bf x}_{ij}\right]^2\\&&=\max_{1\leq i\leq m}\frac{1}{n}\sum_{j=1}^n\left[\varepsilon_{ij}^2-2(\beta-\beta^0)'{\bf x}_{ij}\varepsilon_{ij}+(\beta-\beta^0)'{\bf x}_{ij}{\bf x}_{ij}'(\beta-\beta^0)\right].\end{eqnarray*}
Note that ${\bf x}_{ij},i=1,\cdots,m,j=1,\cdots,n$, are independent and identically distributed with zero mean and variance $\sigma_i^2$. By comparing the asymptotic orders of every terms in the above expression, we have
\begin{eqnarray*}&&\max_{1\leq i\leq m}\frac{1}{n}\sum_{j=1}^n\left[Y_{ij}-\beta'{\bf x}_{ij}\right]^2\\&&=\frac{1}{n}\sum_{j=1}^n\left[\varepsilon_{i_*j}^2-2(\beta-\beta^0)'{\bf x}_{i_*j}\varepsilon_{i_*j}+(\beta-\beta^0)'{\bf x}_{i_*j}{\bf x}_{i_*j}'(\beta-\beta^0)\right]+\delta_n.\end{eqnarray*} As was shown that $\varepsilon_{i_*j}$ and $\delta_n$ are independent of $\beta$. Thus, to get the estimator of $\beta$, minimizing $\max\limits_{1\leq i\leq m}\frac{1}{n}\sum_{j=1}^n\left[Y_{ij}-\beta'{\bf x}_{ij}\right]^2$ is equivalent to minimizing
$$\sum_{j=1}^n\left[-2(\beta-\beta^0)'{\bf x}_{i_*j}\varepsilon_{i_*j}+(\beta-\beta^0)'{\bf x}_{i_*j}{\bf x}_{i_*j}'(\beta-\beta^0)\right].$$ We rewrite the above objective function as
$$Z_n(\gamma)=\sum_{j=1}^n\left[-2\frac{\varepsilon_{ij}{\bf x}_{i_*j}'}{\sqrt{n}}\gamma+\gamma'\frac{{\bf x}_{i_*j}{\bf x}_{i_*j}'}{n}\gamma\right].$$
The function $Z_n(\gamma)$ is obviously convex and is minimized at
$\hat\gamma_n=\sqrt{n}(\hat\beta_G-\beta^0).$ It follows from the Lindeberg-Feller central limit theorem and {\it C1} that
$$Z_n(\gamma)\stackrel{d}\longrightarrow Z_0(\gamma)=-2W'\gamma+\gamma'E({\bf x}{\bf x}')\gamma,$$ where $W\sim N(0,\sigma^2_{i_*}E({\bf x}{\bf x}'))$. The convexity of the limiting objective function, $Z_0(\gamma)$, assures the uniqueness of the minimizer and consequently, that
$$\sqrt{n}(\hat\beta_G-\beta)=\hat\gamma_n=\arg\min\tilde Z_n(\gamma)\stackrel{d}\longrightarrow\hat\gamma_0=\arg\min Z_0(\gamma).$$
(See, e.g., Pollard 1991, Hj\o rt and Pollard 1993, Knight 1998). Finally, we see $\hat\gamma_0=(E({\bf x}{\bf x}'))^{-1}W$ and the result follows. $\Box$

{\bf Proof of Theorem 3.2} The definitions of the two estimations lead directly to the conclusions of the theorem. $\Box$

{\bf Proof of Theorem 3.3} From the proof of Theorem 3.1 we see that $\hat\beta_G$ is actually the common LS estimator of $\beta$ obtained by data $(Y_{i_*j},{\bf x}_{i_*j}),j=1,\cdots,n$. Thus $\hat\beta_G=\beta+O_p(1/n)$, where $\beta$ is the true regression coefficient given by (\ref{(2.5)}) in the mean-certainty model. When
$E[{\bf x}]=0$, the true regression coefficients in the mean-certainty model and the mean-uncertainty model are the same as given in (\ref{(2.5)}) and (\ref{(2.8)}). Moreover, by the the same argument as used in the proof of Theorem 3.1, we have
$$\max_{1\leq i\leq m}\frac 1n\sum_{j=1}^n\left[Y_{ij}-\beta'{\bf x}_{ij}\right]=\max_{1\leq i\leq m}E_{F_i}[\varepsilon]+O_p(1/n).$$
The above discussion ensures that
$$ \hat{\overline\mu}=\max_{1\leq i\leq m}\frac 1n\sum_{j=1}^n\left[Y_{ij}-\beta'{\bf x}_{ij}\right]+O_p(1/n)=\max_{1\leq i\leq m}E_{F_i}[\varepsilon]+O_p(1/n)=\overline\mu+O_p(1/n),$$ where $F_i$ is the distribution of data in $I_i$. Consequently,
$$\max_{1\leq i \leq m}\frac 1n\sum_{j=1}^n\left[Y_{ij}-\beta'{\bf x}_{ij}-\hat{\overline\mu}\right]^2=\max_{1\leq i \leq m}\frac 1n\sum_{j=1}^n\left[Y_{ij}-\beta'{\bf x}_{ij}-\overline\mu\right]^2+O_p(1/n).
$$ On the other hand,
\begin{eqnarray*}&&\frac{1}{n}\sum_{j=1}^n\left[Y_{ij}-\beta'{\bf x}_{ij}-\overline\mu\right]^2=\frac{1}{n}
\sum_{j=1}^n\left[\varepsilon_{ij}-\mu_{i}-(\overline\mu-\mu_{i})\right]^2
\\&&=\sigma_{i}^2+(\overline\mu-\mu_i)^2+O_p(1/n).\end{eqnarray*}
Then,
$$\max_{1\leq i \leq m}\frac{1}{n}\sum_{j=1}^n\left[Y_{ij}-\beta'{\bf x}_{ij}-\overline\mu\right]^2=\frac{1}{n}\sum_{j=1}^n\left[Y_{k_*j}-\beta'{\bf x}_{k_*j}-\overline\mu\right]^2+O_p(1/n).$$
By the above result, $E[{\bf x}]=0$ and the same argument as used in the proof of Theorem 3.1, we can prove the theorem.
$\Box$

{\bf Proof of Theorem 3.4} The proof of the theorem follows directly from the definitions of the two estimators. $\Box$

\

\leftline{\large\bf References}

\begin{description}

\item Artzner, Ph., Delbaen, F., Eber, J. M. and Heath, D. (1997). Thinking
coherently. {\it RISK}, {\bf 10}, 86-71.

\item Bradley, R. C. and Bryc, W. (1985). Multilinear forms and
measures of dependence between random variables. {\it J.
Multivariate Anal.}, {\bf 16}, 335-367.

\item Briand, Ph., Coquet, F., Hu, Y., M\'emin J. and Peng, S. (2000). A converse
comparison theorem for BSDEs and related properties of $g$-expectations.
{\it Electron. Comm. Probab}, {\bf 5}, 101-117.

\item Cand\'es, E. J. and Tao, T. (2007). The Dantzig selector:
statistical estimation when $p$ is much larger than $n$. {\it Ann.
Statist.} {\bf 35}, 2313-2351.

\item Chen, Z. and Epstein, L. (2002). Ambiguity, risk and asset returns in continuous
time. {\it Econometrica}, {\bf 70}(4), 1403-1443.

\item Chen, Z. and Peng, S. (2000). A general downcrossing inequality for gmartingales.
{\it Statist. Probab. Lett.}, {\bf 46}(2), 169-175.

\item Coquet, F., Hu, Y., M\'emin J. and Peng, S. (2002). Filtration-consistent
nonlinear expectations and related $g$-expectations. {\it Probab. Theory Relat.
Fields}, {\bf 123}, 1-27.

\item Denis, L. and Martini, C. (2006). A theoretical framework for the pricing of
contingent claims in the presence of model uncertainty. {\it The Ann. of Appl.
Probability}, {\bf 16}(2), 827-852.

\item Denis, L., Hu, M. and Peng S. (2011). Function spaces and capacity related to a sublinear expectation: application to $G$-Brownian motion pathes. {\it Potential Anal}., {\bf 34}, 139-161.

\item El Karoui, N., Peng, S. and Quenez, M.C. (1997). Backward stochastic
differential equation in finance. {\it Mathematical Finance}, {\bf 7}(1): 1-71.

\item Fan, J. and Li, R. (2001). Variable selection via nonconcave penalized likelihood and its oracle properties. {\it J Amer Statist Assoc}, {\bf 96}, 1348-1360.

\item Fan, J. and Peng, H. (2004). Nonconcave penalized likelihood
with a diverging number of parameters. {\it Ann. Statist.}, {\bf
32}, 928-961.

\item F\"ollmer, H. and Schied, A. (2004). {\it Statistic Finance, An introduction in
discrete time} (2nd Edition), Walter de Gruyter.

\item Gao, F. Q. (2009). Pathwise properties and homeomorphic flows for stochastic
differential equations driven by $G$-Brownianmotion. {\it Stochastic Processes
and their Applications}, {\bf 119}, 3356-3382.

\item Huber,P. J. (1981). {\it Robust Statistics}, John Wiley \& Sons.

\item Hj\o rt, N. and D. Pollard (1993), Asymptotics for minimizers of convex processes.
{\it Statistical Research Report}.

\item Kolmogorov, A. N. and Rozanov, U. A. (1960). On the strong mixing conditions of a
stationary Gaussian process. {\it Probab. Theory Appl.}, {\bf 2},
222-227.

\item Knight, K. (1989). Limit theory for autoregressive-parameter estimates in an infinite-variance random walk. {\it Canadian Journal of Statistics}, {\bf 17}, 261-278.

\item Li, X and Peng, S. (2011). Stopping times and related It\^o's calculus with $G$-Brownian motion. {\it Stochastic Processes and their Applications}, {\bf 121}, 1492-1508.

\item
Lu, C. R. and Lin, Z. Y. (1997). Limit theories for mixing dependent
variables. Science Press, Beijing.

\item Peng, S. (1997). Backward SDE and related $g$-expectations, in Backward
Stochastic Differential Equations, Pitman Research Notes in Math. Series,
No.364, El Karoui Mazliak edit. 141-159.

\item Peng, S. (1999). Monotonic limit theorem of BSDE and nonlinear decomposition
theorem of Doob-Meyer¡¯s type. {\it Prob. Theory Rel. Fields}, {\bf 113}(4),
473-499.

\item Peng, S. (2004). Filtration consistent nonlinear expectations and evaluations
of contingent claims. {\it Acta Mathematicae Applicatae Sinica}. English Series
{\bf 20}(2), 1-24.

\item Peng, S. (2005). Nonlinear expectations and nonlinear Markov chains, {\it Chin.
Ann. Math.}, {\bf 26}B(2), 159-184.

\item Peng, S. (2006). $G$-Expectation, $G$-Brownian Motion and Related
Stochastic Calculus of It\^o¡¯s type, The Abel Symposium 2005, Abel Symposia
2, Edit. Benth et. al., 541-567, Springer-Verlag, 2006.


\item Peng, S. (2008). Multi-dimensional $G$-Brownian motion and related
stochastic calculus under G-expectation. {\it Stochastic Processes and their
Applications}, {\bf 118}(12), 2223-2253.

\item Peng, S. (2009). Survey on normal distributions, central limit theorem,
Brownian motion and the related stochastic calculus under sublinear expectations.
{\it Science in China Series} A: Mathematics, {\bf 52},
7, 1391-1411.


\item Pollard, D. (1991). Asymptotics for least absolute deviation regression Estimators. {\it
Econometric Theory}, {\bf 7}, 186-199.

\item Rosazza, G. E. (2006). Some examples of risk measures via $g$-expectations.
{\it Insurance: Mathematics and Economics}, {\bf 39}, 19-34.

\item Rosenblatt, M. (1956). A central limit theorem and a strong
mixing condition. {\it Proc. Nat. Acad. Sci.}, {\bf 42}, 43-47.

\item Rosenblatt, M. (1970). Density estimates and Markov sequences,
in {\it Nonparametric Techniques in Statistical inference}, ed. M.
Puri, London: Cambridge University Press, pp. 199-210.

\item Soner, M., Touzi, N. and Zhang, J. (2011a). Martingale representation theorem under G-expectation. {\it Stochastic Processes and their Applications}, {\bf 121},
    265-287.

\item  Soner M, Touzi N, Zhang J. (2011b) Quasi-sure stochastic analysis through aggregation. {\it Electronic Journal of Probability}, {\bf 16}, 1844-1879.

\item  Soner M, Touzi N, Zhang J. (2012). Wellposedness of second order backward SDEs. {\it Probability Theory and Related Fields}, {\bf 153}, 149-190.

\item  Soner M, Touzi N, Zhang J. (2013). Dual formulation of second order target problems. {\it Annals of Applied Probability}, {\bf 23}, 308-347.

\item Song, Y.(2012). Uniqueness of the representation for $G$-martingales with finite variation. {\it Electron. J.Probab}., {\bf 17}, 1-15.

\item Strassen, V. (1964). Messfehler und information, {\it Z. Wahrsheinlichkeitstheorie verw. Gebiete}, {\bf 2}, 267-284.

\item Tibshirani, R. J. (1996). Regression shrinkage and selection via the Lasso. {\it J. R. Stat. Soc. Ser. B}, {\bf 58}, 267-288.

\item Walley, P. (1991). {\it Statistical Reasoning with Imprecise Probabilities},
Chapman and Hall, London, New York.

\item Xu, J. and Zhang, B. (2009). Martingale characterization of $G$-Brownian motion.
{\it Stochastic Processes and their Applications}, {\bf 119}, 232-248.

\item Zou, H. (2006). The adaptive lasso and its oracle properties. {\it J. Amer. Statist. Assoc.,} {\bf 101}, 1148-1429.

\end{description}

\newpage

%

\end{document}